\documentclass[english]{article}
\usepackage[T1]{fontenc}
\usepackage[latin9]{inputenc}
\usepackage{geometry}
\usepackage{verbatim}
\usepackage{float}
\usepackage{amsmath}
\usepackage{graphicx}
\usepackage{esint}

\makeatletter

\newcommand{\lyxdot}{.}

\newcommand{\lyxaddress}[1]{
\par {\raggedright #1
\vspace{1.4em}
\noindent\par}
}

\makeatother

\usepackage{babel}
\begin{document}

\title{Condition Numbers of Indefinite Rank $2$ Ghost Wishart Matrices}

\author{Ramis Movassagh%
\thanks{ramis@math.mit.edu%
} $\;$and Alan Edelman%
\thanks{edelman@math.mit.edu%
} }

\maketitle

\lyxaddress{\begin{center}
Department of Mathematics, Massachusetts Institute of Technology,
Cambridge MA, 02139
\par\end{center}}
\begin{abstract}
We define an indefinite Wishart matrix as a matrix of the form $A=W^{T}W\Sigma$,
where $\Sigma$ is an indefinite diagonal matrix and $W$ is a matrix
of independent standard normals. We focus on the case where $W$ is
$L\times2$ which has engineering applications. We obtain the distribution
of the ratio of the eigenvalues of $A$. This distribution can be
``folded'' to give the distribution of the condition number. We
calculate formulas for $W$ real $\left(\beta=1\right)$, complex
$\left(\beta=2\right)$, quaternionic $\left(\beta=4\right)$ or any
ghost $0<\beta<\infty$. We then corroborate our work by comparing
them against numerical experiments. 
\end{abstract}

\section*{Problem Statement}

Let $W$ be an $L\times2$ matrix whose elements are drawn from a
normal distribution. Let the two real eigenvalues of $A=W^{T}W\Sigma$
be denoted by $\lambda_{1}$ and $\lambda_{2}$. The condition number
of $A$ is $\sigma=\frac{\left|\underset{}{\lambda_{\mbox{max}}}\right|}{\left|\lambda_{\mbox{min}}\right|}$,
where $\left|\lambda_{\mbox{max}}\right|=\mbox{max}\left(\left|\lambda_{1}\right|,\left|\lambda_{2}\right|\right)$
and $\left|\lambda_{\mbox{min}}\right|=\mbox{min}\left(\left|\lambda_{1}\right|,\left|\lambda_{2}\right|\right)$.
What is the condition number distribution of $A=W^{T}W\Sigma$, where
$\Sigma=\left[\begin{array}{cc}
x_{1}\\
 & x_{2}
\end{array}\right]$ is full rank with $\mbox{sgn}\left(x_{2}\right)=-\mbox{sgn}\left(x_{1}\right)$?
Though much is known when $\Sigma$ is positive definite \cite{muirhead},
to the extent of our knowledge, the indefinite case is rarely considered.
We work out the distribution of the ratio of the eigenvalues and the
condition number of $A$ as it has applications in hypersensitive
ground based radars \cite{Christ,Christ2}. 

A ghost $\beta-$normal is a construction that works in many ways
like a real $\left(\beta=1\right)$, complex $\left(\beta=2\right)$,
and quaternionic $\left(\beta=4\right)$ standard normal \cite{Alan_Ghosts}.
In particular, its absolute value has a $\chi_{\beta}$ distribution
on $\left[0,\infty\right)$ . In general, we allow $W$ to be an $L\times2$
matrix sampled from a $\beta-$normal distribution but we immediately
turn this into a problem involving real matrices namely the condition
number of $R\Sigma R^{T}$, where $R\sim\left[\begin{array}{cc}
\chi_{L\beta} & \chi_{\beta}\\
 & \chi_{(L-1)\beta}
\end{array}\right]$.

\section*{Distribution of the Ratio of the Eigenvalues }

We first write $W=QR$, where $Q$ is an $L\times2$ orthogonal matrix
and $R$ is an $2\times2$ upper triangular matrix \cite[reduced QR factorization]{Nick}.
Note that $W\Sigma W^{T}=QR\Sigma R^{T}Q^{T}$ and that $A$ is similar
to $R\Sigma R^{T}$. The elements of $R$ may be chosen to be non-negative
in which case it is well known that they have \textit{independent}
$\chi$-distributions 
\begin{equation}
R\equiv\left[\begin{array}{cc}
a & b\\
 & c
\end{array}\right]\sim\left[\begin{array}{cc}
\chi_{L\beta} & \chi_{\beta}\\
 & \chi_{(L-1)\beta}
\end{array}\right]\label{eq:R_chiDistributions}
\end{equation}
 and that \cite{muirhead} 
\begin{equation}
\chi_{k}\sim\frac{x^{k-1}e^{-x^{2}/2}}{2^{k/2-1}\Gamma(k/2)},\quad x\ge0\label{eq:chi}
\end{equation}
where $k$ denotes degrees of freedom and $\beta=1$ corresponds to
entries being real, $\beta=2$ complex, $\beta=4$ quaternionic and
general $\beta$ for any ghost \cite{Alan_Ghosts}. The positivity
of $R_{1,2}$ merits some comment as it does not generalize well beyond
two columns. With two columns, even for $\beta\ne1$, a phase can
be pulled out of the rows and columns of $R$ and absorbed elsewhere
without any loss of generality. 

Comment: The concept of a $QR$ decomposition is often sensible for
$W$'s with entries drawn from $\beta-$normal distribution. As with
reals, complex, and quaternions, $R$ can be chosen to have positive
diagonal elements for rank $2$ matrices, and $R_{12}$ can be chosen
non-negative as well, by absorbing \textquotedbl{}phases\textquotedbl{}
or unit-ghosts either into $Q$ or $R\Sigma R^{T}$.

This argument follows the ghost methodology described in \cite{Alan_Ghosts}.
However, mathematically the point of departure can be the computation
of the condition number distribution for the real matrix $R\Sigma R^{T}$
without any mention of the ghosts or their interpretation.

The joint distribution of the elements $a,b,c$ of $R$ with parameters
$L$ and $\beta$ is:

\begin{equation}
\rho_{R}\left(a,b,c;L,\beta\right)=\frac{2^{-\beta L+3}a^{\beta L-1}b^{\beta-1}c^{\beta\left(L-1\right)-1}\exp\left[-\frac{1}{2}\left(a^{2}+b^{2}+c^{2}\right)\right]}{\Gamma\left(\frac{\beta L}{2}\right)\Gamma\left(\frac{\beta}{2}\right)\Gamma\left[\frac{\beta}{2}\left(L-1\right)\right]};\label{eq:Rdensity}
\end{equation}
in particular 

\begin{eqnarray*}
\rho_{R}\left(a,b,c;L,\beta=1\right) & = & \frac{2a^{L-1}c^{L-2}\exp\left[-\frac{1}{2}\left(a^{2}+b^{2}+c^{2}\right)\right]}{\pi\Gamma\left(L-1\right)},\\
\rho_{R}\left(a,b,c;L,\beta=2\right) & = & \frac{2^{-2L+3}\left(L-1\right)a^{2L-1}bc^{2L-3}\exp\left[-\frac{1}{2}\left(a^{2}+b^{2}+c^{2}\right)\right]}{\Gamma^{2}\left(L\right)},\\
\rho_{R}\left(a,b,c;L,\beta=4\right) & = & \frac{2^{-4\left(L-1\right)}\left(L-1\right)\left(2L-1\right)a^{4L-1}b^{3}c^{4L-5}\exp\left[-\frac{1}{2}\left(a^{2}+b^{2}+c^{2}\right)\right]}{\Gamma^{2}\left(2L\right)}\mbox{ }.
\end{eqnarray*}

The first change of variables computes the matrix whose condition
number we are seeking, $R\Sigma R^{T}$ as: $R\left[\begin{array}{cc}
x_{1}\\
 & x_{2}
\end{array}\right]R^{T}=\left[\begin{array}{cc}
a^{2}x_{1}+b^{2}x_{2} & bcx_{2}\\
bcx_{2} & c^{2}x_{2}
\end{array}\right]\equiv\left[\begin{array}{cc}
d & f\\
f & e
\end{array}\right]$. From Eq. (\ref{eq:R_chiDistributions}) we see that $a,b,c$ are
real and non-negative. The old and new variables are related by 

\begin{equation}
c=\sqrt{e/x_{2}},\qquad b=\frac{f}{x_{2}\sqrt{e/x_{2}}},\qquad a=\sqrt{\frac{1}{x_{1}}\left(d-\frac{f^{2}}{e}\right)}.\label{eq:abs}
\end{equation}
Note that $c>0$ implies $\mbox{sgn}\left(e\right)=\mbox{sgn}\left(x_{2}\right)$
and $f=bcx_{2}\Rightarrow\mbox{sgn}\left(f\right)=\mbox{sgn}\left(x_{2}\right)$.
The Jacobian associated with this transformation is $\left|\frac{\partial(a,b,c)}{\partial(d,e,f)}\right|=\left|\frac{1}{4ac^{2}x_{1}x_{2}^{2}}\right|$. 

We make the eigenvalue decomposition of the symmetric matrix

\begin{equation}
\left[\begin{array}{cc}
d & f\\
f & e
\end{array}\right]\equiv\left[\begin{array}{cc}
\cos\theta & -\sin\theta\\
\sin\theta & \cos\theta
\end{array}\right]\left[\begin{array}{cc}
\lambda_{1}\\
 & \lambda_{2}
\end{array}\right]\left[\begin{array}{cc}
\cos\theta & \sin\theta\\
-\sin\theta & \cos\theta
\end{array}\right],\quad\theta\in\left[0,\frac{\pi}{2}\right];
\end{equation}
implying 
\begin{eqnarray}
d & = & \lambda_{1}\cos^{2}\theta+\lambda_{2}\sin^{2}\theta\label{eq:d}\\
e & = & \lambda_{1}\sin^{2}\theta+\lambda_{2}\cos^{2}\theta\label{eq:e}\\
f & = & \frac{\mbox{sin}2\theta}{2}\left(\lambda_{1}-\lambda_{2}\right).\label{eq:f}
\end{eqnarray}
 The Jacobian associated with this transformation is $\left|\frac{\partial\left(d,e,f\right)}{\partial\left(\lambda_{1},\lambda_{2},\theta\right)}\right|=|\lambda_{1}-\lambda_{2}|$.
The choice of $\theta\in\left[0,\frac{\pi}{2}\right]$ nails down
the ordering of the eigenvalues

\begin{equation}
\begin{array}{ccccc}
x_{2}<0<x_{1} & \Rightarrow & f<0 & \Rightarrow & \lambda_{1}<0<\lambda_{2}\\
x_{1}<0<x_{2} & \Rightarrow & f>0 & \Rightarrow & \lambda_{2}<0<\lambda_{1}.
\end{array}\label{eq:lambdas}
\end{equation}
In summary, given $\theta\in\left[0,\frac{\pi}{2}\right]$, the constraints
on $d,e,f$ are

\[
\begin{array}{ccccc}
x_{2}<0<x_{1} & \Rightarrow & d>\frac{f^{2}}{e},\mbox{ }e<0,\mbox{ }f<0, & \Rightarrow & \lambda_{1}<0<\lambda_{2}\\
x_{1}<0<x_{2} & \Rightarrow & d<\frac{f^{2}}{e},\mbox{ }e>0,\mbox{ }f>0 & \Rightarrow & \lambda_{2}<0<\lambda_{1}.
\end{array}
\]

Intersecting $t\equiv-\frac{\lambda_{2}}{\lambda_{1}}\ge0$ with the
constraint space, we obtain 
\begin{equation}
\tan^{2}\theta>t.\label{eq:constraint}
\end{equation}
It is easy to see this is required by Eq. (\ref{eq:e}). Conversely
given $\mbox{sgn}\left(e\right)=\mbox{sgn}\left(x_{2}\right)$ and
$\tan^{2}\theta>t=-\frac{\lambda_{2}}{\lambda_{1}}$ we can solve
for $d$ and $f$ so that $d,e$ and $f$ meet all of the constraints. 

The next change of variables is to write the distribution as a ratio
of the eigenvalues, i.e., $u=\lambda_{1},\: t=\frac{-\lambda_{2}}{\lambda_{1}},$
and $\theta$. The Jacobian associated with this transformation is
$\left|\frac{\partial\left(\lambda_{1},\lambda_{2},\theta\right)}{\partial\left(u,t,\theta\right)}\right|=\left|u\right|$.%
{} Since we are interested in the condition number, we integrate $u$
over $\left\{ \begin{array}{ccc}
\left[0,\infty\right)\; & : & x_{2}>0\\
\left[0,-\infty\right) & : & x_{2}<0
\end{array}\right.\;$, while applying the constraint (\ref{eq:constraint})

\begin{equation}
\rho\left(t,\theta;L,\beta\right)=\left|\frac{2\left(-x_{1}x_{2}\right)^{\frac{\beta L}{2}}\Gamma\left(\beta L\right)t^{\frac{\beta L}{2}-1}\left(t+1\right)^{\beta}}{\Gamma\left(\frac{\beta}{2}\right)\Gamma\left[\frac{\beta}{2}\left(L-1\right)\right]\Gamma\left(\frac{\beta L}{2}\right)}\frac{\left(\sin\theta\cos\theta\right)^{\beta-1}\left[t\cos^{2}\theta-\sin^{2}\theta\right]^{\beta(L-1)}}{\left[x_{1}\left(t^{2}\cos^{2}\theta+\sin^{2}\theta\right)-x_{2}t\right]^{\beta L}}\right|,\label{eq:rho_tbeta_theta}
\end{equation}
where here and below implicitly $x_{1}$ and $x_{2}$ are parameters
as well. This in the special cases reads

\begin{eqnarray*}
\rho\left(t,\theta;L,\beta=1\right) & = & \frac{2\left(-x_{1}x_{2}\right)^{L/2}\Gamma\left(L\right)t^{\frac{L}{2}-1}\left(t+1\right)}{\sqrt{\pi}\Gamma\left(\frac{L}{2}\right)\Gamma\left[\frac{1}{2}\left(L-1\right)\right]}\left|\frac{\left[t\cos^{2}\theta-\sin^{2}\theta\right]^{L-1}}{\left[x_{1}\left(t^{2}\cos^{2}\theta+\sin^{2}\theta\right)-x_{2}t\right]^{L}}\right|\\
\rho\left(t,\theta;L,\beta=2\right) & = & \frac{\left(-x_{1}x_{2}\right)^{L}\Gamma\left(2L\right)t^{L-1}\left(t+1\right)^{2}}{\Gamma\left(L-1\right)\Gamma\left(L\right)}\frac{\mbox{sin}2\theta\left[t\cos^{2}\theta-\sin^{2}\theta\right]^{2(L-1)}}{\left[x_{1}\left(t^{2}\cos^{2}\theta+\sin^{2}\theta\right)-x_{2}t\right]^{2L}}\\
\rho\left(t,\theta;L,\beta=4\right) & = & \frac{\left(-x_{1}x_{2}\right)^{2L}\Gamma\left(4L\right)t^{2L-1}\left(t+1\right)^{4}}{4\Gamma\left[2\left(L-1\right)\right]\Gamma\left(2L\right)}\frac{\mbox{sin}^{3}2\theta\left[t\cos^{2}\theta-\sin^{2}\theta\right]^{4(L-1)}}{\left[x_{1}\left(t^{2}\cos^{2}\theta+\sin^{2}\theta\right)-x_{2}t\right]^{4L}}.
\end{eqnarray*}

Lastly we integrate $\theta$ from $\left[0,\pi/2\right]$ to obtain
the distribution corresponding to the absolute value of the ratio
of the eigenvalues $\rho\left(t;L,\beta\right)$. We denote the density
after integrating $\theta$ with the same symbol $\rho$ without any
confusion. In terms of the Gauss hypergeometric function \cite{TableIntegrals},
$_{2}F_{1}\left(a,b;c;z\right)$ we obtain,

\vspace{0.5cm}

\begin{eqnarray}
\rho(t;L,\beta) & = & \frac{2^{\beta L-1}\left(-x_{1}x_{2}\right)^{\frac{\beta L}{2}}\Gamma\left(\frac{\beta L}{2}+\frac{1}{2}\right)\Gamma\left[\beta\left(L-1\right)+1\right]\left(t+1\right)^{\beta}t^{\frac{\beta L}{2}-\frac{\beta}{2}-1}}{\sqrt{\pi}\Gamma\left[\frac{\beta}{2}\left(L-1\right)\right]\Gamma\left(\beta L-\frac{\beta}{2}+1\right)\left|tx_{2}-x_{1}\right|^{\beta L}}\left(\frac{x_{1}-tx_{2}}{tx_{1}-x_{2}}\right)^{\beta L}\label{eq:rho}\\
 & \times & _{2}F_{1}\left(\beta L,\frac{\beta}{2};\beta L-\frac{\beta}{2}+1;\frac{x_{1}-tx_{2}}{x_{2}-tx_{1}}\right).\nonumber 
\end{eqnarray}

\vspace{0.5cm}
In particular,

\begin{eqnarray*}
\rho\left(t;L,\beta=1\right) & = & \frac{2^{L-1}(-x_{1}x_{2})^{L/2}\Gamma\left(\frac{1}{2}(L+1)\right)\Gamma\left(L\right)\left(t+1\right)t^{\frac{L}{2}-\frac{3}{2}}}{\sqrt{\pi}\Gamma\left(\frac{1}{2}\left(L-1\right)\right)\Gamma\left(L+\frac{1}{2}\right)\left|tx_{2}-x_{1}\right|^{L}}\left(\frac{x_{1}-tx_{2}}{tx_{1}-x_{2}}\right)^{L}\\
 & \times & _{2}F_{1}\left(L,\frac{1}{2};L+\frac{1}{2};\frac{x_{1}-tx_{2}}{x_{2}-tx_{1}}\right),\\
\rho\left(t;L,\beta=2\right) & = & \frac{2^{2\left(L-1\right)}\left(x_{1}t-x_{2}\right)\left(-x_{1}x_{2}\right)^{L}\Gamma\left(L-\frac{1}{2}\right)\left(t+1\right)t^{L-2}}{\sqrt{\pi}\Gamma\left(L-1\right)\left(x_{1}-x_{2}\right)\left(tx_{2}-x_{1}\right)^{2L}}\left(\frac{x_{1}-tx_{2}}{tx_{1}-x_{2}}\right)^{2L},\\
\rho\left(t;L,\beta=4\right) & = & \frac{2^{4L-1}(-x_{1}x_{2})^{2L}\Gamma\left(2L+\frac{1}{2}\right)\Gamma\left(4L-3\right)\left(t+1\right)^{4}t^{2L-3}}{\sqrt{\pi}\Gamma\left[2\left(L-1\right)\right]\Gamma\left(4L-1\right)\left(tx_{2}-x_{1}\right)^{4L}}\left(\frac{x_{1}-tx_{2}}{tx_{1}-x_{2}}\right)^{4L}\\
 & \times & _{2}F_{1}\left(4L,2;4L-1;\frac{x_{1}-tx_{2}}{x_{2}-tx_{1}}\right).
\end{eqnarray*}

\section*{Distribution of the Condition Number}

We can turn $\rho\left(t;L,\beta\right)$ to the condition number
distribution by ``folding'' the distribution about one, i.e., take
$\sigma\equiv\max\left(t,\frac{1}{t}\right)$. Therefore by folding
the answer about one we can form the true condition number distribution
from $\rho(t;L,\beta)$. Mathematically,

\begin{flushleft}
\begin{eqnarray*}
1=\int_{0}^{+\infty}dt\rho(t;L,\beta) & = & \int_{0}^{1}dt\rho\left(t;L,\beta\right)+\int_{1}^{\infty}dt\rho\left(t;L,\beta\right)\overset{\sigma\rightarrow\frac{1}{t}}{=}\int_{1}^{\infty}d\sigma\left(\frac{1}{\sigma^{2}}\rho\left(\frac{1}{\sigma};L,\beta\right)+\rho(\sigma;L,\beta)\right)\\
 & \equiv & \int_{1}^{\infty}d\sigma f(\sigma;L,\beta)
\end{eqnarray*}

\par\end{flushleft}

This way we arrive at the desired distribution function, $f\left(\sigma,L,\beta\right)$,
for the condition number of the indefinite Wishart matrices

\begin{equation}
f\left(\sigma;L,\beta\right)\equiv\frac{1}{\sigma^{2}}\rho\left(\frac{1}{\sigma};L,\beta\right)+\rho\left(\sigma;L,\beta\right),\label{eq:condition}
\end{equation}
where $\rho$ is as in Eq. (\ref{eq:rho}).

\vspace{0.5cm}

The condition number distribution written explicitly is

\begin{eqnarray}
f\left(\sigma;L,\beta\right) & = & \frac{2^{\beta L-1}\left(-x_{1}x_{2}\right)^{\frac{\beta L}{2}}\Gamma\left(\frac{\beta L}{2}+\frac{1}{2}\right)\Gamma\left[\beta\left(L-1\right)+1\right]\left(\sigma+1\right)^{\beta}\sigma^{\frac{\beta L}{2}-\frac{\beta}{2}-1}}{\sqrt{\pi}\Gamma\left[\frac{\beta}{2}\left(L-1\right)\right]\Gamma\left(\beta L-\frac{\beta}{2}+1\right)}\label{eq:general}\\
 & \times & \left[{\scriptstyle \left|\frac{1}{\sigma x_{1}-x_{2}}\right|^{\beta L}{}_{2}F_{1}\left(\beta L,\frac{\beta}{2};\beta L-\frac{\beta}{2}+1;r\right)+\left|\frac{1}{\sigma x_{2}-x_{1}}\right|^{\beta L}{}_{2}F_{1}\left(\beta L,\frac{\beta}{2};\beta L-\frac{\beta}{2}+1;\frac{1}{r}\right)}\right]\nonumber \\
\mbox{where} &  & r\equiv\frac{\sigma x_{2}-x1}{\sigma x_{1}-x_{2}},\qquad\sigma=\left[1,\infty\right).\nonumber 
\end{eqnarray}
\vspace{0.5cm}
In particular, 

\begin{eqnarray}
f\left(\sigma;L,\beta=1\right) & = & \frac{2^{L-1}(-x_{1}x_{2})^{L/2}\Gamma\left(\frac{L}{2}+\frac{1}{2}\right)\Gamma\left(L\right)(\sigma+1)\sigma^{\frac{1}{2}L-\frac{3}{2}}}{\sqrt{\pi}\Gamma\left(\frac{L}{2}-\frac{1}{2}\right)\Gamma\left(\frac{1}{2}+L\right)}\label{eq:real}\\
 & \times & \left[{\scriptstyle \left|\frac{1}{x_{2}-\sigma x_{1}}\right|^{L}{}_{2}F_{1}\left(\frac{1}{2},L;L+\frac{1}{2};r\right)+\left|\frac{1}{\sigma x_{2}-x_{1}}\right|^{L}{}_{2}F_{1}\left(\frac{1}{2},L;L+\frac{1}{2};\frac{1}{r}\right)}\right],\nonumber \\
f\left(\sigma;L,\beta=2\right) & = & \frac{1}{4}\frac{(-4x_{1}x_{2})^{L}\Gamma(L-\frac{1}{2})(\sigma+1)}{\sqrt{\pi}(x_{1}-x_{2})\Gamma(L-1)\sigma^{4}}\left[\frac{\sigma^{L+2}(\sigma x_{1}-x_{2})}{(x_{2}-x_{1}\sigma)^{2L}}+\frac{\sigma^{L+2}(x_{1}-\sigma x_{2})}{(x_{2}\sigma-x_{1})^{2L}}\right],\label{eq:complex}\\
f\left(\sigma;L,\beta=4\right) & = & \frac{2^{4L-1}\left(-x_{1}x_{2}\right)^{2L}\Gamma\left(2L+\frac{1}{2}\right)\Gamma\left[4\left(L-1\right)+1\right]\left(\sigma+1\right)^{4}\sigma^{2L-3}}{\sqrt{\pi}\Gamma\left[2\left(L-1\right)\right]\Gamma\left(4L-1\right)}\label{eq:quaternion}\\
 & \times & \left[{\scriptstyle \left(\frac{1}{\sigma x_{1}-x_{2}}\right)^{4L}{}_{2}F_{1}\left(4L,2;4L-1;r\right)+\left(\frac{1}{\sigma x_{2}-x_{1}}\right)^{4L}{}_{2}F_{1}\left(4L,2;4L-1;\frac{1}{r}\right)}\right].\nonumber 
\end{eqnarray}

\section*{Special Cases }

Consider the case where $x_{1}=-x_{2}\equiv x$ which implies $r=-1$
and the eigenvalues of $A$ and $-A$ are equi-distributed . The condition
number density is 

\[
f\left(\sigma;L,\beta\right)=\frac{2\sigma^{\frac{\beta L}{2}-\frac{\beta}{2}-1}\Gamma\left(\beta L-\beta+1\right)\left(\sigma+1\right)^{-\beta\left(L-1\right)}}{\beta\left(L-1\right)\left\{ \Gamma\left[\frac{\beta}{2}\left(L-1\right)\right]\right\} ^{2}}
\]
and in particular takes the very simple form

\begin{equation}
f\left(\sigma;L=2,\beta=1\right)=\frac{2}{\pi\sqrt{\sigma}(\sigma+1)}.\label{eq:simpleCase}
\end{equation}

Alternatively as $\beta\rightarrow\infty$ the limiting condition
number is non-random and equal to the condition number of the matrix

\[
\left[\begin{array}{cc}
\sqrt{L} & 1\\
 & \sqrt{L-1}
\end{array}\right]\left[\begin{array}{cc}
x_{1}\\
 & x_{2}
\end{array}\right]\left[\begin{array}{cc}
\sqrt{L} & 1\\
 & \sqrt{L-1}
\end{array}\right]^{T}.
\]
This is a simple consequence of the fact that $\lim_{s\rightarrow\infty}\frac{1}{\sqrt{s}}\chi_{s}=1$
and the condition number is independent of scaling.

\section*{Numerical Results}

In this section we compare the theoretical condition number distribution,
Eq. (\ref{eq:general}), (red curves) against Monte Carlo data (black
dots) for the three special cases of real, complex and quaternionic
matrices $W$. In all the plots the number of trials used to generate
Monte Carlo data was $10^{5}$.

\subsection*{Real Matrices $\beta=1$}

\begin{figure}[H]
\includegraphics[scale=0.45]{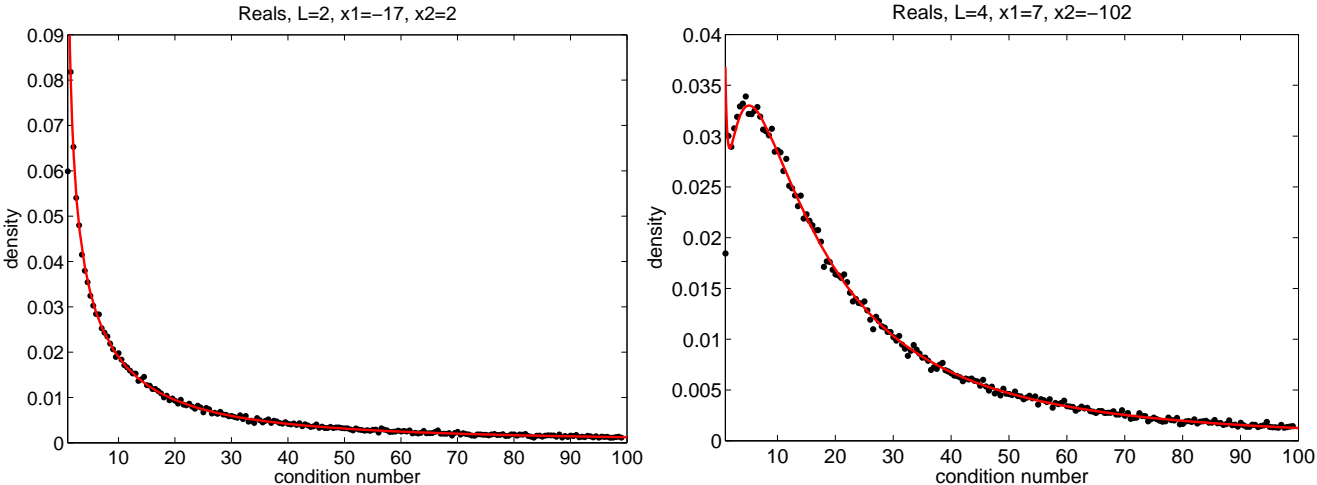}\includegraphics[scale=0.47]{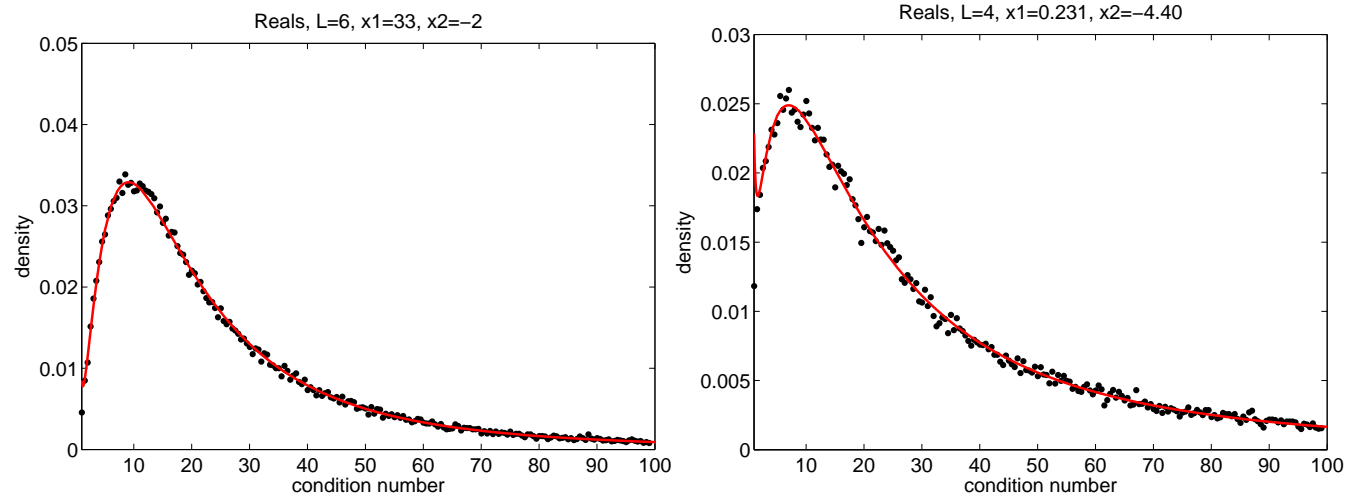}

\caption{Real $W$ matrices, $\beta=1$, given by Eq. (\ref{eq:real})}
\end{figure}

\begin{figure}[H]
\includegraphics[scale=0.47]{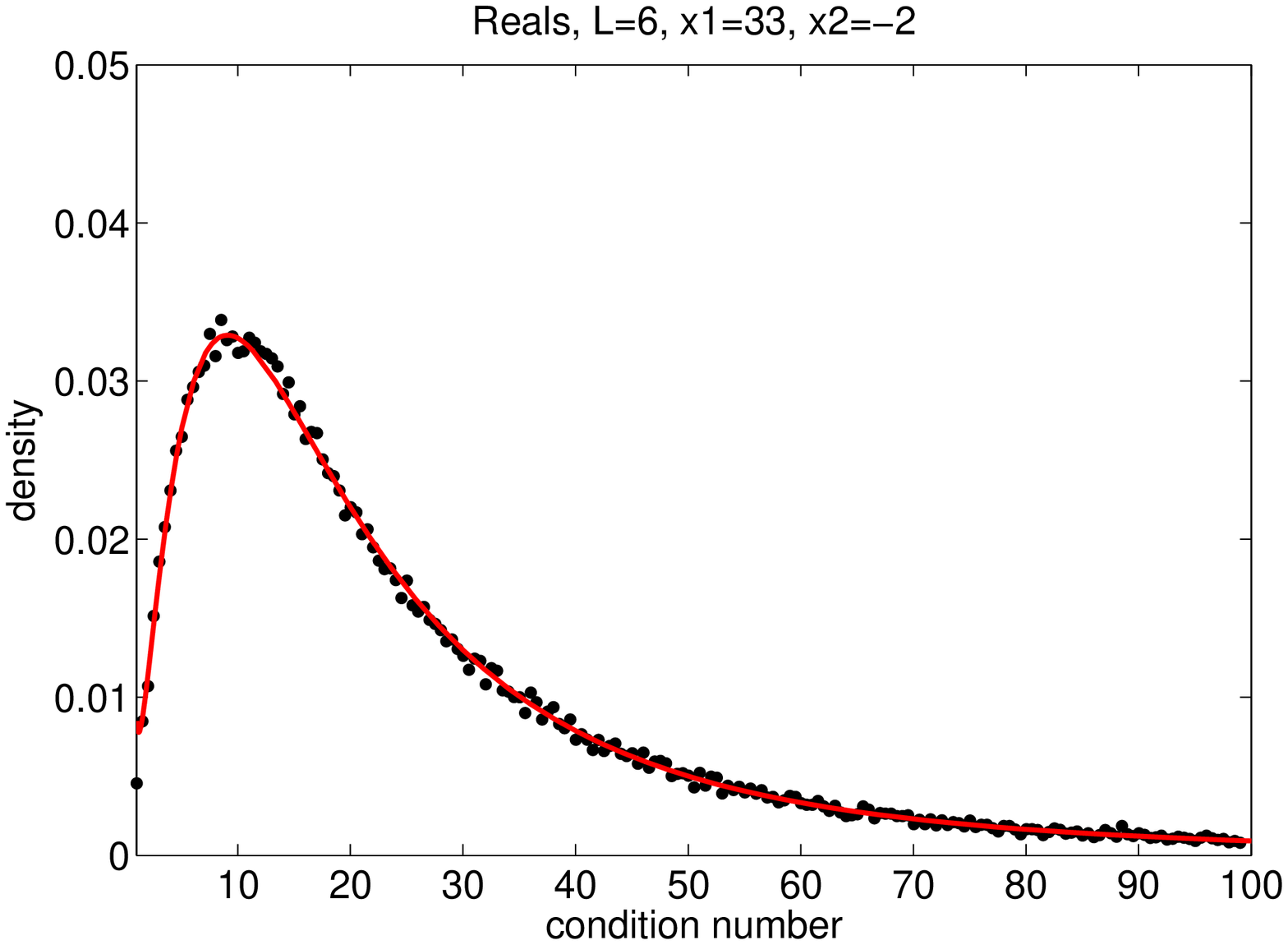}\includegraphics[scale=0.47]{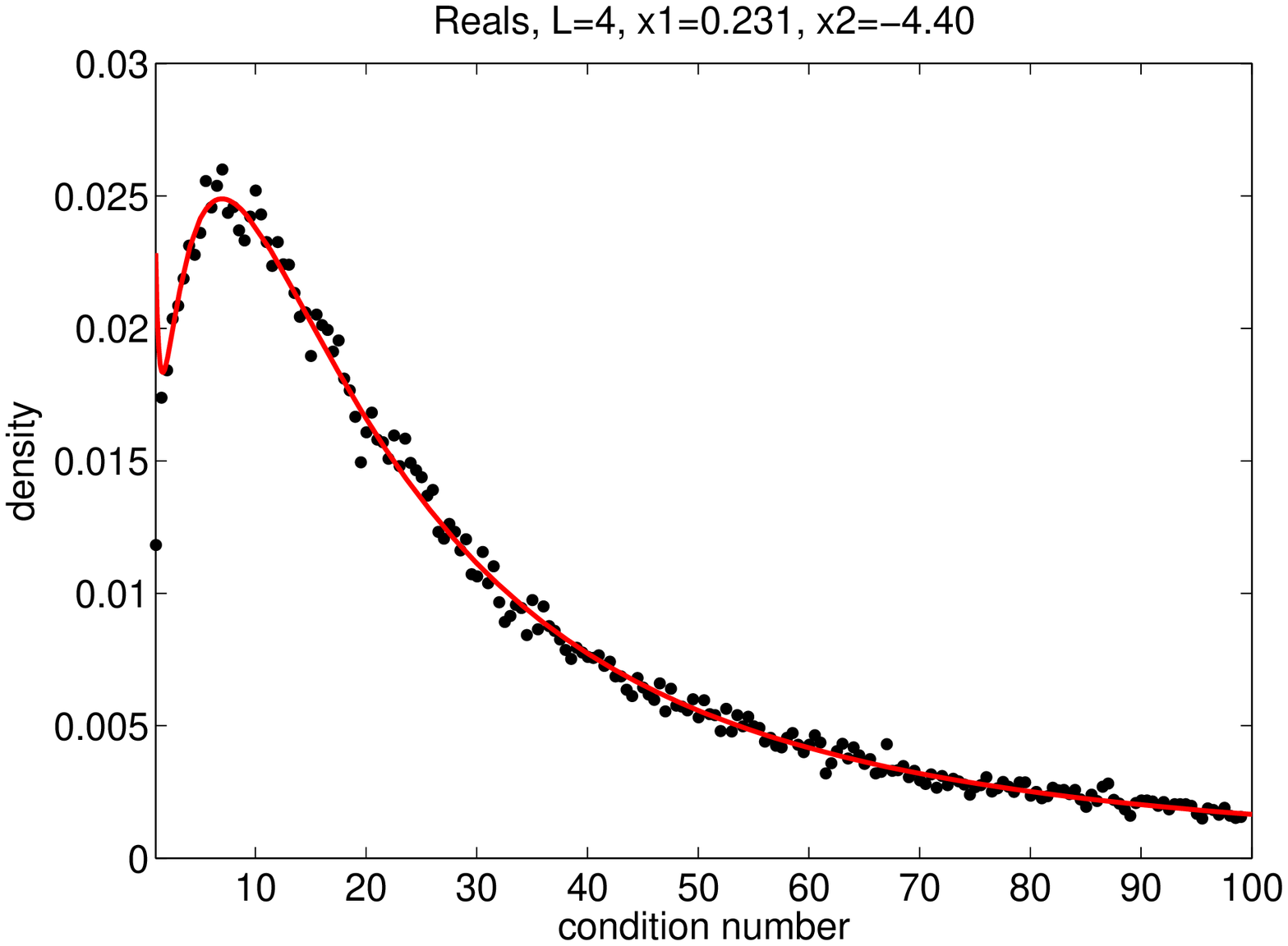}\caption{Real $W$ matrices, $\beta=1$, given by Eq. (\ref{eq:real})}

\end{figure}

\subsection*{Complex Matrices $\beta=2$}

\begin{figure}[H]
\includegraphics[scale=0.47]{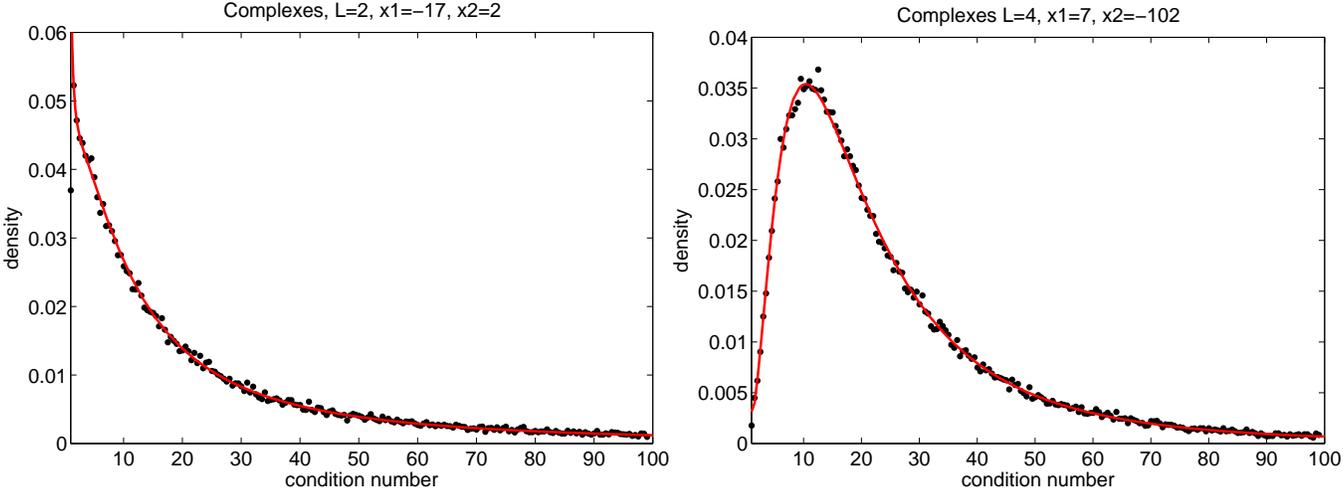}\includegraphics[scale=0.47]{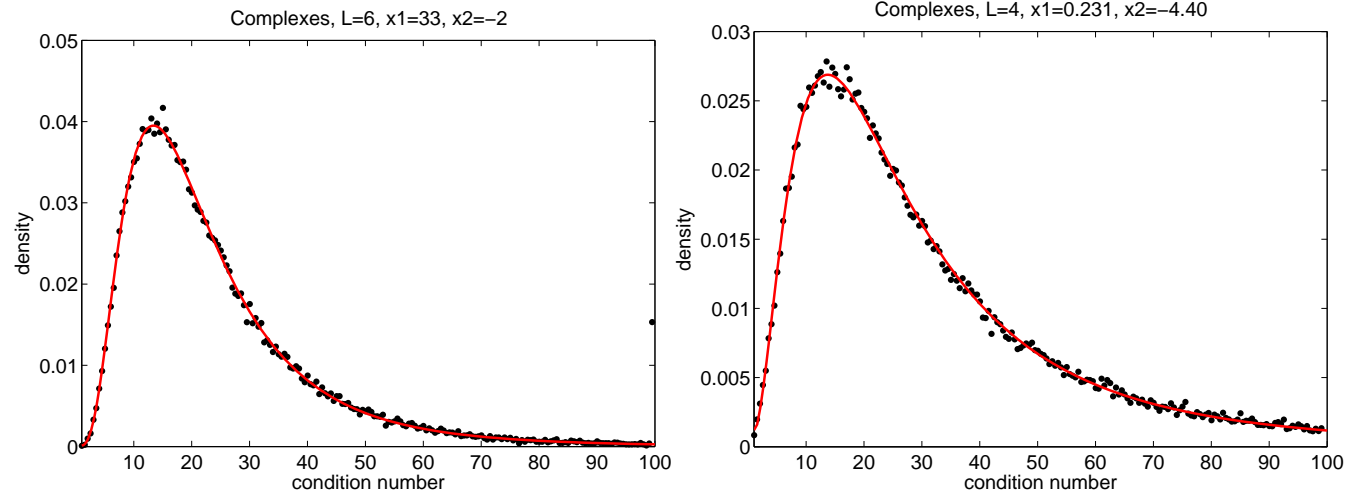}

\caption{Complex $W$ matrices, $\beta=2$, given by Eq. (\ref{eq:complex})}

\end{figure}

\begin{figure}[H]
\includegraphics[scale=0.47]{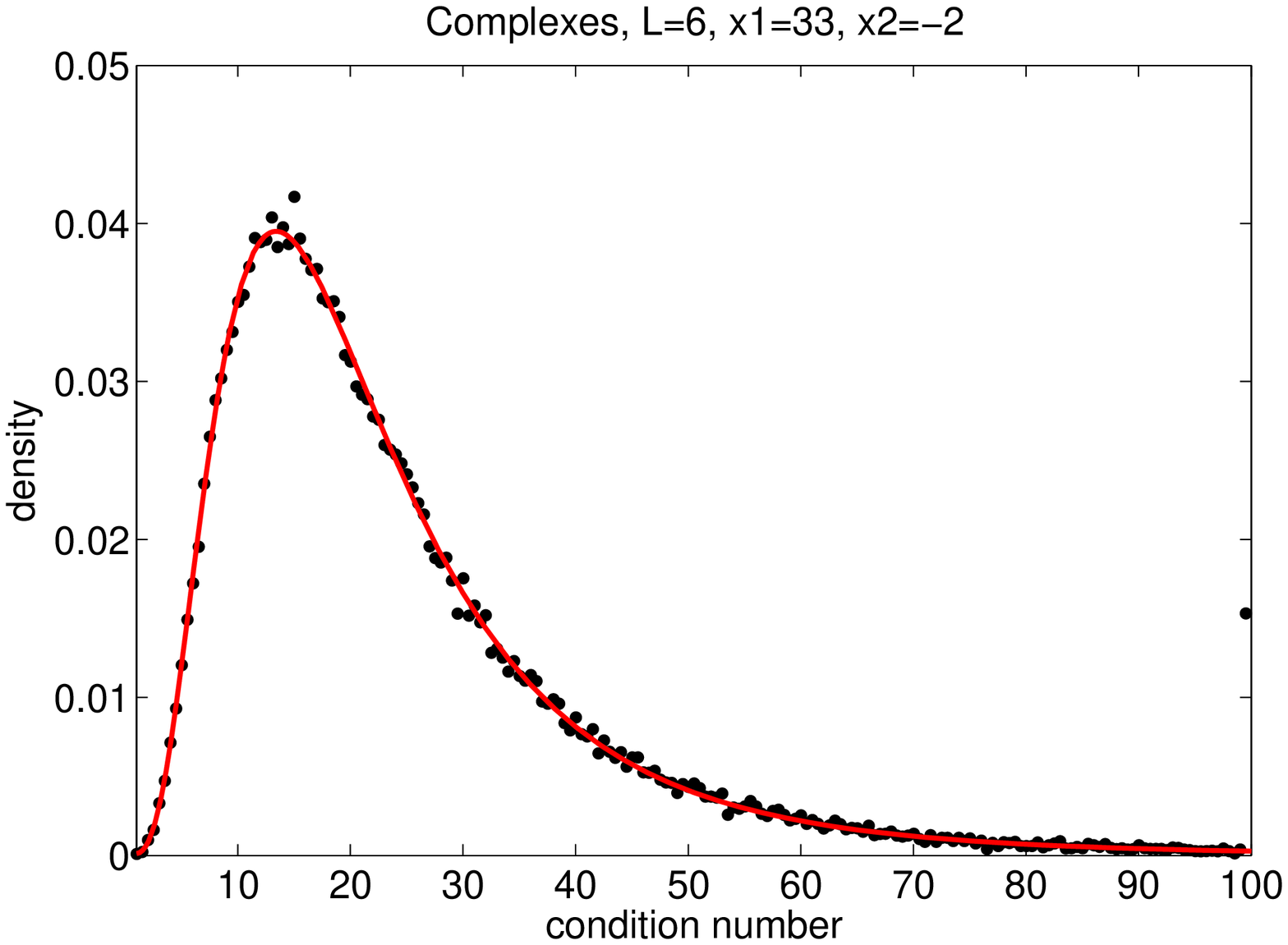}\includegraphics[scale=0.47]{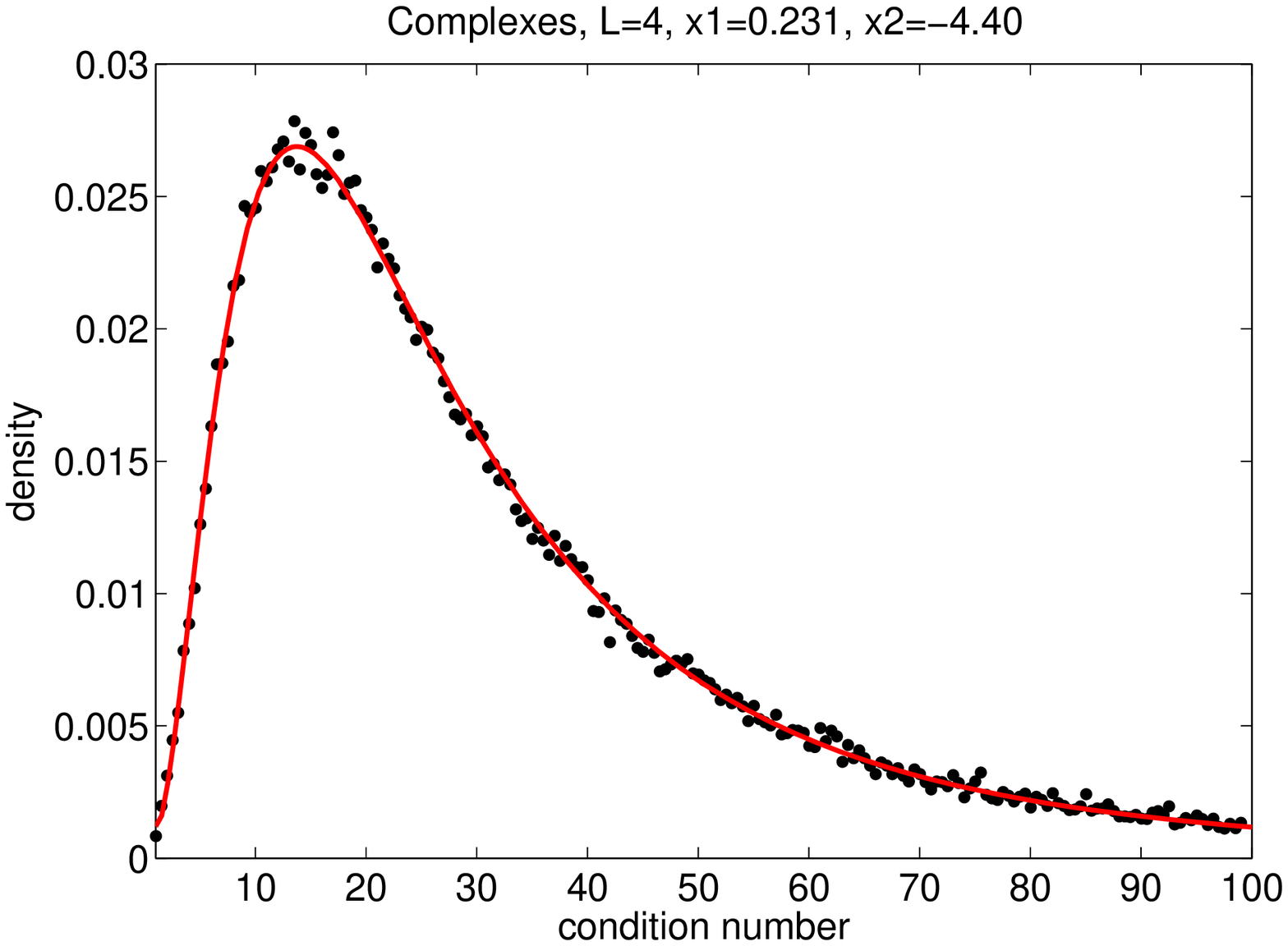}\caption{Complex $W$ matrices, $\beta=2$, given by Eq. (\ref{eq:complex})}

\end{figure}

\subsection*{Quaternionic Matrices $\beta=4$}

\begin{figure}[H]
\includegraphics[scale=0.47]{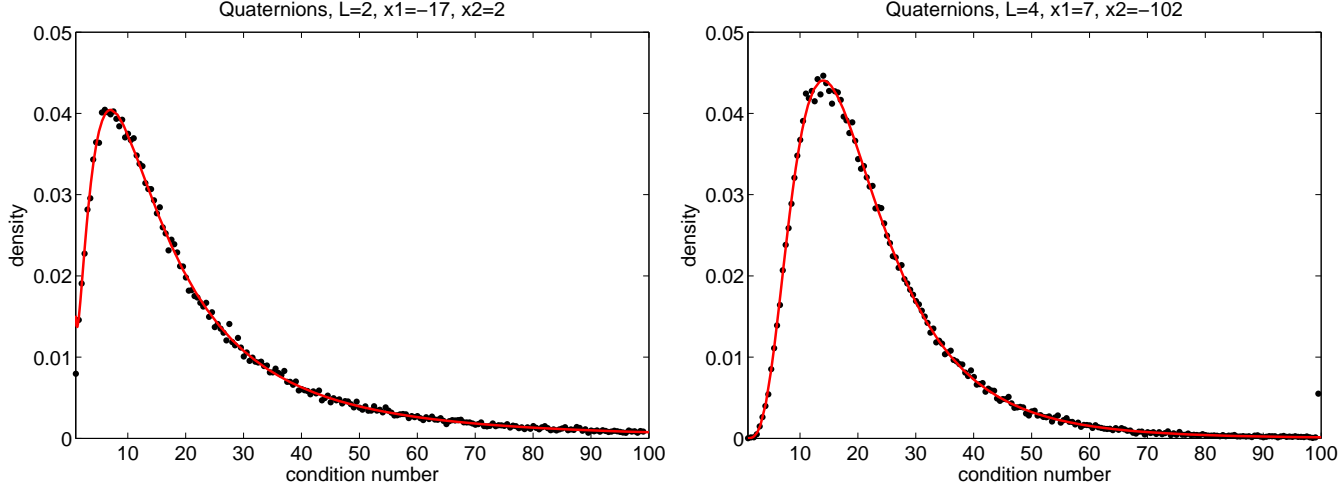}\includegraphics[scale=0.47]{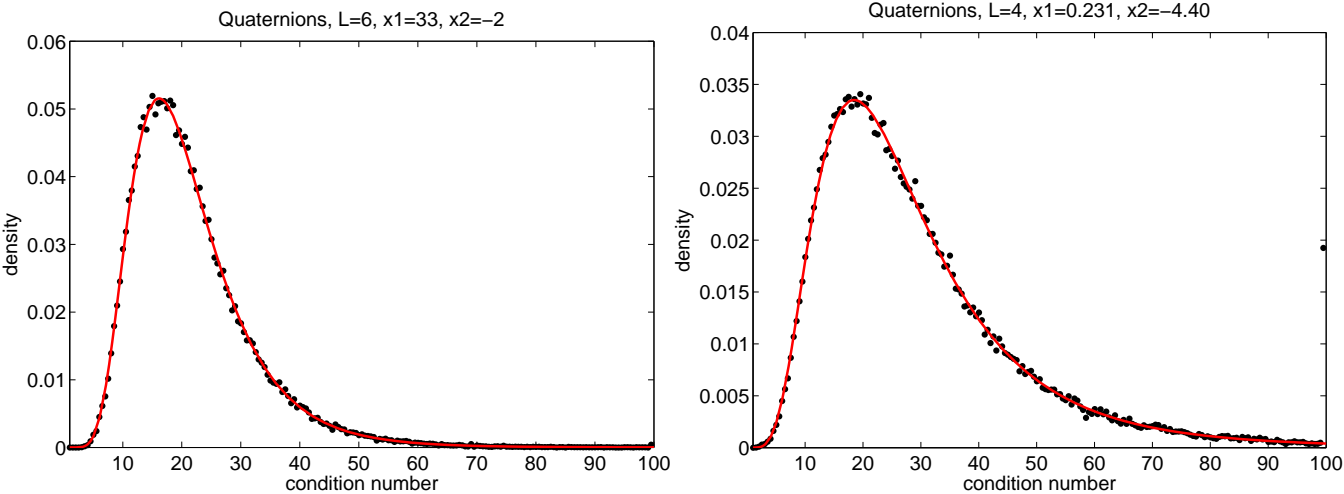}

\caption{Quaternionic $W$ matrices, $\beta=4$, given by Eq. (\ref{eq:quaternion})}
\end{figure}

\begin{figure}[H]
\includegraphics[scale=0.47]{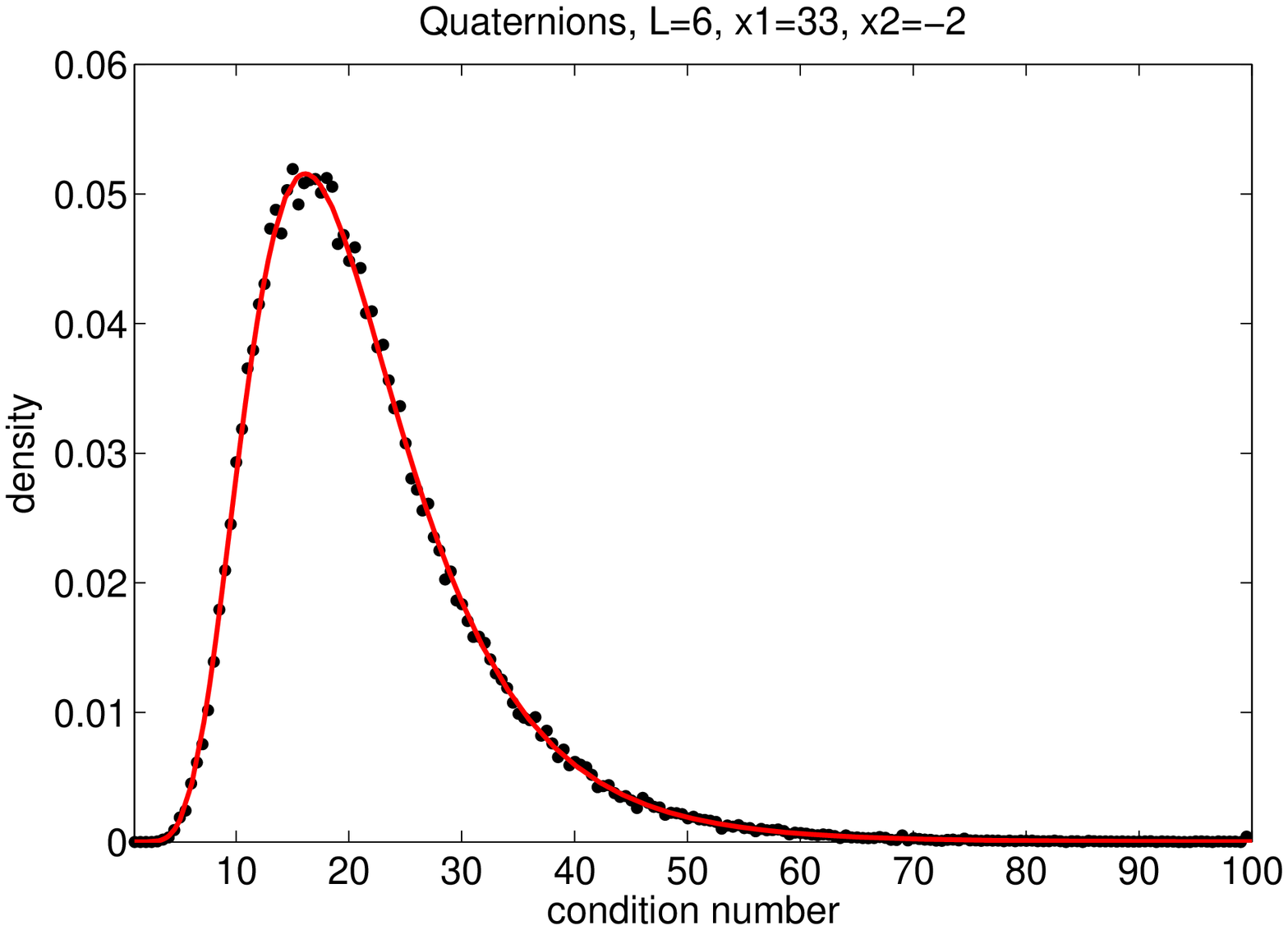}\includegraphics[scale=0.47]{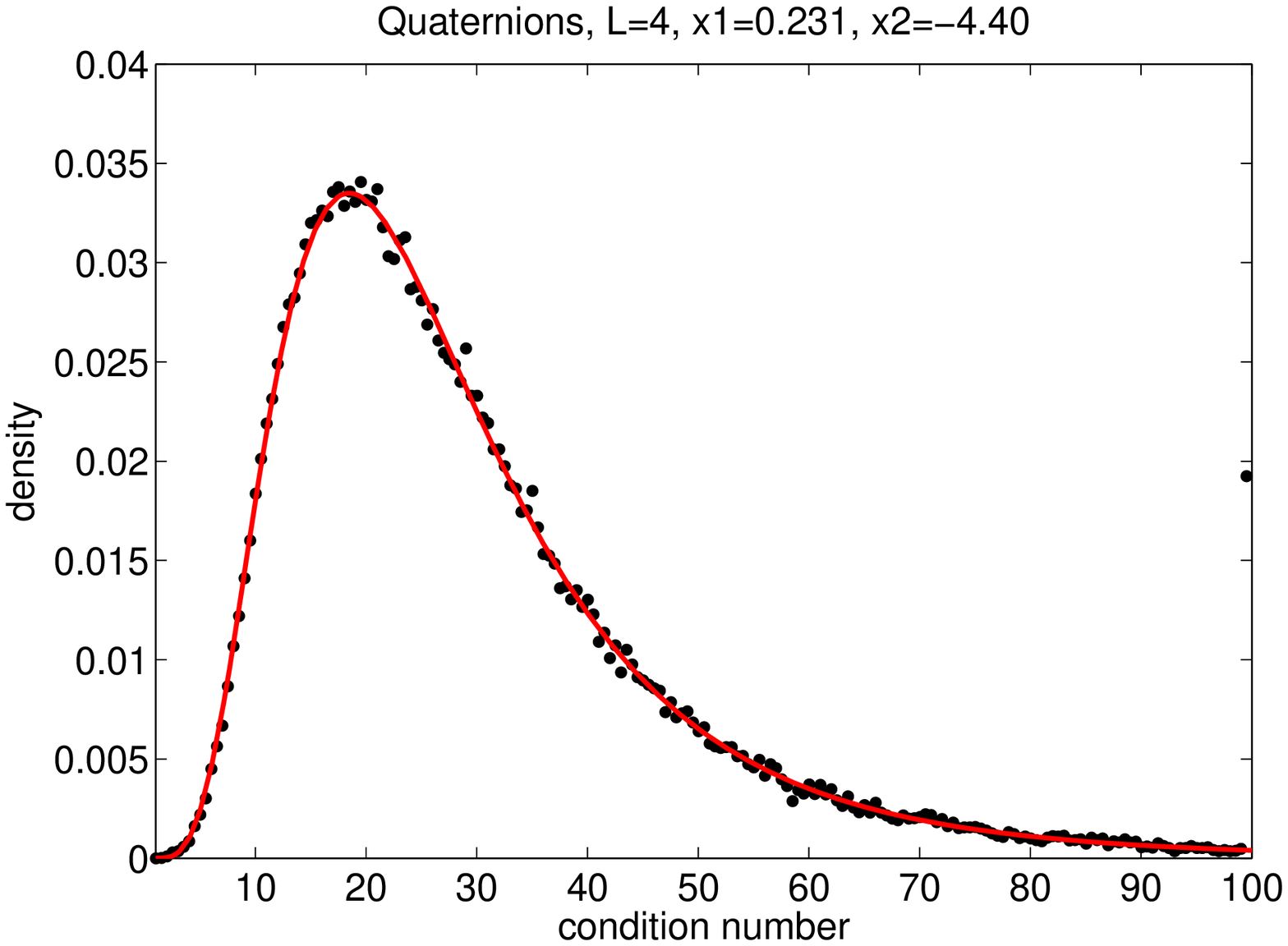}\caption{Quaternionic $W$ matrices, $\beta=4$, given by Eq. (\ref{eq:quaternion})}
\end{figure}

\subsection*{General $\beta$}

Here we plot Eq. (\ref{eq:general}) for a fixed set of parameters
$L,x_{1},x_{2}$ for various ghosts.

\begin{figure}[H]
\begin{centering}
\includegraphics[scale=0.6]{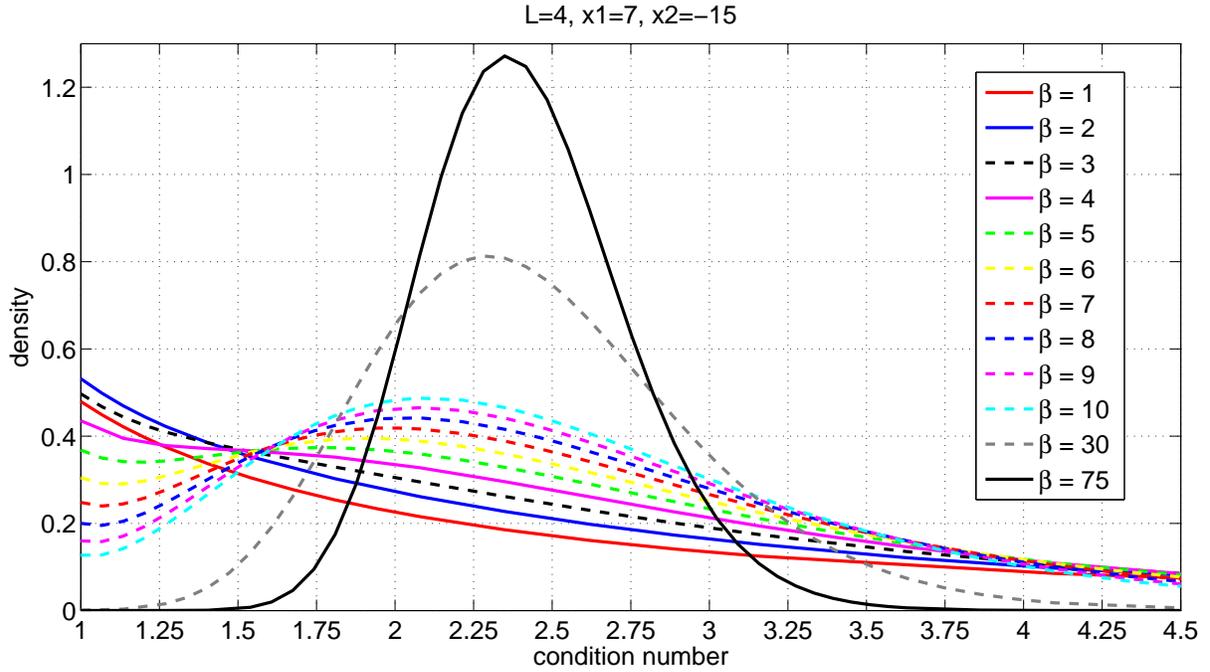}\caption{The condition number density Eq. (\ref{eq:general}) for various ghosts.
As $\beta\rightarrow\infty$ the distribution becomes non-random and
can be represented by a delta function around $2.395$.}

\par\end{centering}

\end{figure}

\section*{Future Work}

The natural generalization would be to extend the results to $L\times N$
indefinite matrices. We satisfied ourselves with $N=2$ given its
relevance for applications \cite{Christ,Christ2}.

\section*{Acknowledgements}

We thank Christ Richmond, Alexander Dubbs, and Peter W. Shor for discussions.
Maple and Mathematica were used for symbolic integrations and MATLAB
for the generation of Monte Carlo data and plotting. We acknowledge
the support of National Science Foundation through grants  CCF-0829421,
SOLAR Grant No. 1035400, and DMS-1035400.

\end{document}